\documentclass[a4,13pt]{amsart}
\usepackage{amssymb,amstext,amsmath,amscd,amsthm,amsfonts,enumerate,graphicx,latexsym, hyperref,cleveref}
\usepackage[usenames]{color}
\usepackage[all]{xy}

 \usepackage{comment}
 \usepackage{mathrsfs}
 \usepackage{mathtools}
 
 \RequirePackage{doi}
\usepackage{hyperref}

\newtheorem{thm}{Theorem}[section]

\newtheorem{cor}[thm]{Corollary}

\newtheorem{lem}[thm]{Lemma}
\newtheorem{prop}[thm]{Proposition}
\newtheorem{ques}[thm]{Question}

\newtheorem*{theorem*}{Theorem}

\theoremstyle{definition}
\newtheorem{dfn}[thm]{Definition}
\newtheorem{rem}[thm]{Remark}

\newtheorem{ex}[thm]{Example}

\newtheorem*{claim*}{Claim}
\theoremstyle{remark}
\newtheorem{setup}[thm]{Setup}

\numberwithin{equation}{thm}

\def\cond{\mathfrak c}

\def\End{\operatorname{End}}

\def\Hom{\operatorname{Hom}}

\def\im{\operatorname{Im}}

\def\m{\mathfrak{m}}

\def\Tor{\operatorname{Tor}}

\def\tr{\operatorname{tr}}

\def\ad{\operatorname{add}}
\newcommand{\floor}[1]{\lfloor #1 \rfloor}

\DeclarePairedDelimiter{\brackets}{(}{)}

\newcommand{\En}[2]{\operatorname{End}_{#1} (#2)}

\renewcommand{\hom}[3]{\operatorname{Hom}_{#1} (#2, #3)}

\newcommand{\Z}{\operatorname{{Z}} \brackets}

\begin{document}
\setlength{\baselineskip}{15pt}
\title{Stable trace ideals and applications}
\author{Hailong Dao}
\address{Department of Mathematics, University of Kansas, Lawrence, KS 66045-7523, USA}
\email{hdao@ku.edu}
\urladdr{https://www.math.ku.edu/~hdao/}

\author{Haydee Lindo}
\address{Department of Mathematics, Harvey Mudd College, 301 Platt Blvd.
Claremont, CA 91711, USA}
\email{hlindo@hmc.edu}

\thanks{The authors are supported by the National Science Foundation Grant No. 2137949 and  Simons Collaborator Grant FND0077558.}


\begin{abstract}
We study stable trace ideals in  one dimensional local Cohen-Macaulay rings and give numerous applications. 
\end{abstract}
\maketitle

\section{Introduction}
Let $R$ be a one dimensional Cohen-Macaulay local ring. Let $I$ be an ideal of height one in $R$. Recall that  $I$  is called {\it stable} if $I\cong \End_R(I)$ as $R$-modules. $I$ is called a {\it trace ideal} if $I= \sum_{f\in \Hom_R(I,R)} \im(f)$. Stable ideals can be thought of as ideals with simplest blow-ups, and Lipman exploited their nice properties in his  seminal  (\cite{lipman1971stable}) on Arf rings. Trace ideals have long been  useful technical tools in commutative algebra, but recently they have attracted new attention as interesting objects in their own right, see \cite{FaberTrace, haydee_trace, GotoIsobeKumashiro, herzog2019trace}. 

In this paper, we study stable trace ideals in detail. Our motivation comes from many sources. First, via endomorphism rings, they are in a one-to-one correspondence with the set of finite birational extensions of $R$ that are self-dual as $R$-modules. These birational extensions play a key role in studying reflexive modules over $R$, see \cite{daoArf}. If the integral closure of $R$ is finite, then its conductor is always a stable trace ideal. Second, a key result in \cite{haydee_trace} on modules whose endomorphism ring has Gorenstein center can be viewed as a statement about modules whose trace ideal is stable, and our collaboration started from this insight. 

We soon realize that stable trace ideals satisfy even more delightful properties in general, and they seem to hold the key to understand reflexive modules, trace ideals and integrally closed ideals in general. We hope that our work will serve as a starting point for more studies in these directions. 

We now describe the key results and organization of our paper. After a preliminary section, we focus on basic property of stable and stable trace ideals in \Cref{sec_st}. We shall focus on regular ideals, those that contain a regular element. We prove that if either $I$ or the trace of $I$, $\tr I$ is stable, then $\tr I$ is isomorphic to $I^*=\Hom_R(I,R)$, the $R$-dual of $I$.

In \Cref{sec_rf}, we point out the connections between stable ideals and reflexive ideals. For instance, if $\tr I$ is stable, then $I$ is reflexive if and only $I\cong I^*$ if and only if $I\cong \tr I$. We also give a characterization for when $\tr I\cong I^*$ using the conductor of the birational  blow-up ring of $I$.

Next we establish two key theorems in \Cref{sec_2thm}. The first one says that if $M$ is a module who trace $I=\tr M$ is stable and regular, then $I$ is a direct summand of a direct sum of the dual $M^*$. Thus, if Krull-Schmidt holds, then $I$ is a direct summand of $M^*$. If furthermore $M$ is reflexive, then $I$ must be a summand of $M$. The second theorem says that if $I,J$ and $I\cap J$ are stable trace ideals, then $\Hom_R(I,J)\cong\Hom_R(J,I) \cong I\cap J \cong \tr{(IJ)}$. This result generalizes the well-known fact that if $I$ is a stable trace ideal, then $\Hom_R(I,I)\cong I \cong I^2\cong I^*$.

In \Cref{sec_int} we study the connections between stable trace ideals and integrally closed ideals. Fix a regular integrally closed ideal $J$ and consider the set $T(J)$ of ideals $I$ whose integral closure is $J$. The key result here states that among $T(J)$, the stable ones must be minimal with respect to inclusion. The concept of $T(J)$ seems to be interesting in its own right, and it clearly plays in important role in understand the abundance (or lack thereof) of trace ideals.  We are able to compute it in several examples. 

The last two sections give some brief applications of the theory in some currently active topics. For instance,  in \Cref{SectionArf} we give a new characterization of Arf rings: they are the ones where any regular trace ideal is stable. Furthermore, over a complete reduced Arf ring, any reflexive modules with a rank decompose into direct sum of stable trace ideals. Finally, \Cref{sec_gor} focuses on classes  rings defined by the trace of canonical ideals. For instance, we  give an alternative proof of a result in \cite{herzog2019trace} that when $R$ has minimal multiplicity, the properties of being almost Gorenstein and nearly Gorenstein are equivalent. We also recover a result in \cite{herzog2021tiny} that characterizes canonical ideals whose trace equals the conductor.

\section{Preliminaries}
In this section we recall some basic notions and results needed for subsequent sections. 

Throughout $R$ is a commutative Noetherian ring. We write $Q(R)$ for the total quotient ring of $R$. Given an $R$-module $M$ and an ideal $I$ we write $M^*$ for $\hom R M R$, $I^{-1}$ for the set $\{x \in Q(R) \, | \, xI \subseteq R \}$ and $\bar{I}$ for the integral closure of $I$.

\begin{dfn}
The \emph{trace map} on $M$ is the homomorphism $\tau_M: M \otimes M^* \longrightarrow R$ given by $m \otimes \alpha \mapsto \alpha(m)$. The \emph{trace ideal of $M$},  denoted $\tr(M)$,  is the image of $\tau_M$ in $R$. We say an ideal $I$ of $R$ is \emph{trace} (or is a \emph{trace ideal}) provided there exists an $R$-module $M$ such that $I =\tr(M)$. 
\end{dfn}

\begin{rem}
Note that if $M\cong N$ as $R$-modules then $\tr{M }= \tr{N}.$
\end{rem}

\begin{lem}
If $I$ is a trace ideal then the following hold true
\begin{enumerate}
    \item $\tr(I) = I$
    \item The inclusion $I\subseteq R$ induces the identification $I^* = \En R I$. 
\end{enumerate}
\end{lem}
\begin{proof}
c.f. \cite[Proposition 2.8]{haydee_trace}. 
\end{proof}

The following lemma is well-known; see, for example, \cite[Proposition 4.14]{WLCMR} for a similar result to Lemma \ref{homRS} in slightly different language. We provide proofs for the sake of completion. 

\begin{lem}\label{homRS} 
Let $S$ be a ring such that $R\subseteq S \subset Q(R)$. Suppose $M$ and $N$ are $S$-modules. If $N$ is torsionfree $S$-module, then \[\hom R M N = \hom S M N . \]
\end{lem}

\begin{proof}
Since $R\subseteq S$ all $S$-linear homomorphisms are $R$-linear and $\hom S M N  \subseteq \hom R M N $.  Given any $\alpha \in \hom R M N$, $m \in $M, and  $\dfrac{a}{b} \in S $ note that \[ b\left[\dfrac{a}{b} \alpha (m) - \alpha \left( \dfrac{a}{b} m \right) \right] = 0 .\]

Since $b$ is a nonzerdivisor and $N$ is torsionfree we conclude that  $\dfrac{a}{b} \alpha(m)  - \alpha \left( \dfrac{a}{b} m \right) = 0$ and therefore $\alpha$ is $S$-linear.\end{proof}

\section{Stable trace ideals in dimension one}\label{sec_st}

Throughout this section we assume $(R,\m)$ is a local  Cohen-Macaulay ring of dimension one. Recall that an ideal $I$ is \emph{regular} if it contains a nonzerodivisor. The blow-up ring of $I$ in $R$ is defined as $B_R(I):= \cup_{n>0} I^n:_{Q(R)}I^n$. Such a regular ideal $I$ is called \emph{stable} if  $B_R(I)= I:_{Q(R)}I$ (\cite[Definition 1.3]{lipman1971stable}).

In what follows we write  $\Z{\En R I}$ for the center of the endomorphism ring of an ideal $I$. For a regular ideal $I$, we identify the modules  $I^*: = \hom R I R$ and $I^{-1}:= \{x \in Q(R) \,| \, xI \subseteq R \}$  as subsets of $Q(R)$.

\begin{lem} \label{isoend}
Let $I$ and $J$ be regular ideals in $R$. If $J \cong I$ then $\End_R(I) = \End_R(J)$ in $Q(R)$.
\end{lem}

\begin{proof}
Because these ideals contains a nonzerodivisor, $J \cong I$ implies $J = xI$ for some $x \in Q(R)$. For any $y \in \End_R(I)$, \[yJ = yxI = xyI \subseteq xI = J.\] We have shown that $\End_R(I) \subseteq \End_R(J)$. A symmetric argument proves the other containment. 
\end{proof}

The following is largely known, but for lack of convenient references we provide  proofs. 

\begin{prop}\label{prop_stable}
Let $I$ be a regular ideal. The following are equivalent: 
\begin{enumerate}
\item $I$ is stable.
\item $I\cong \End_R(I)$.
\item $I^2=xI$ for some regular $x\in I$.
\item $I\cong I^2$.
\end{enumerate}
\end{prop}

\begin{proof}
Let $S= I:_{Q(R)}I$. Clearly $S\cong \End_R(I)$. Since $I$ is an ideal of $S$ and $S$-isomorphisms are $R$-isomorphisms between the two (\Cref{homRS}), $(2)$ is equivalent to $I=xS$ for some $x\in S$. Note that such $x$ must be in $I$ and regular.  If $I=xS$, then $I^2=x^2S^2=x^2S=xI$, so $(2)$ implies $(3)$. Assume $(3)$, then $\frac{I^2}{x}=I$, thus $\frac{I}{x}\subset S$. But clearly $S\subset \frac{I}{x}$, hence $S=\frac{I}{x}\cong I$, and we established the equivalence of $(2)$ and $(3)$. Finally $(3)$ implies $(4)$, and if $(4)$ holds then $I^n\cong I$ for any $n$, thus $B_R(I)= S$ by \Cref{isoend}. The equivalence of $(1)$ and $(3)$ is \cite[Lemma 1.11]{lipman1971stable}.
\end{proof}

The next result collects the main characterizations of stable trace ideals. 

\begin{prop} \label{prop_trace}

Let $I$ be a regular trace ideal. The following are equivalent: 
\begin{enumerate}
\item $I^2=xI$ for some regular element $x\in I$.
\item $I^2\subset (x)$ for some regular element $x\in I$. 
\item $I=(x):I$ for some regular element $x\in I$. 
\item $I\cong I^{*}$.
\item $I\cong \End_R(I)$. \label{stabletraceend}
\item $I\cong I^2$.  
\item $I$ is stable.
\end{enumerate}


\end{prop}

\begin{proof}
The equivalence of $(1), (5), (6), (7)$ is \Cref{prop_stable}.
$(1)\implies (2)$ is clear. 

$(2) \implies (3)$: if $I^2 \subseteq (x)$, note that \[I \subset (x): I = xI^{-1} \subseteq I I^{-1} = \tr(I) = I.\] 
The claim follows.

$(3) \implies (4)$: Note $I =  (x): I = xI^{-1} = x I^* \cong I^*$.

$(4) \implies (5)$: This follows from the assumption that $I$ is regular and trace which implies an equality in $Q(R)$ between $I^*$ and $\End_R(I)$.

\end{proof}

\begin{rem}
If $I$ in Proposition \ref{prop_trace} has a principal reduction, by \cite[Proposition 4.5]{dms} we may take $(x)$ in the proposition to be that reduction. The assumption that an ideal has a principal reduction is a mild one. For example, it is true for all ideals whenever $R$ has an infinite residue field; see \cite[ 8.3.7, 5.1.6]{SwansonHuneke}.
\end{rem}

\begin{lem}\label{lem_red}
If $A\subseteq J$ is a reduction of a regular ideal $J$ then $AJ^*$ is a reduction of $\tr(J)$ and there exists an $n$ such that $A\tr(J)^n= J\tr(J)^n$. 
\end{lem}

\begin{proof}

Note, $AJ^* \subseteq JJ^* = \tr(J).$
As a reduction, there exists some $n$ such that $AJ^n = J^{n+1}$. Then in $Q(R)$ \[ AJ^* \tr(J)^n = AJ^* J^n (J^*)^{n} = J^{n+1}(J^*)^{n+1} =\tr(J)^{n+1}. \] 

Also, 
\[A \tr(J)^n = AJ^n (J^*)^{n}= J^{n+1}(J^*)^{n}= J \tr(J)^n \qedhere \]
\end{proof}

\begin{rem}
By this lemma, when $J$ has a principal reduction, we have $J^* \tr(J)^n \cong \tr(J)^{n+1} $ for some $n$.  Proposition \ref{IvJ} below is concerned with the case $n=0$. It considers results of the kind obtained in Proposition \ref{prop_trace} in the case where an ideal $I$ is stable but not necessarily a trace ideal. 
\end{rem}

\begin{lem}\label{traceformula}(\cite[Lemma 3.6]{dms})
Suppose that $I$ is a regular ideal and $x\in I$ be a non zero divisor. Then $\tr(I)= I((x):_R I):_R x$. \qedhere \qed
\end{lem}

\begin{prop} \label{IvJ}
Let $J$ be an  regular ideal with a principal reduction $(x)$ and let $I=\tr(J)$. Consider the following:
\begin{enumerate}
    \item $I\cong I^2$.
    \item $I=(x):J$.
    \item $I\cong (x):J$.
    \item $I\cong J^{*}$.
\end{enumerate}
\end{prop}
Then $(1)\implies (2)\iff (3)\iff (4)$.
\begin{proof}

Clearly (3) and (4) are equivalent and (2) implies both.  Assume (4). Let $J'=(x):J \cong J^*$. Then by \ref{traceformula} $JJ' = xI \cong I\cong J'$. Since $J$ has a principal reduction and $J'$ is MCM (since it is a dual over a one-dimensional ring), by \cite[Prop 4.5]{dms} and \Cref{traceformula}, $xI= JJ' = xJ'$, thus $I=J'$.

$(1) \implies (4)$: by \Cref{lem_red} there is $n$ such that $xI^n= JI^n$. Since $I$ is stable, it follows that $\frac{J}{x}\subset I^n:_{Q(R)}I^n = I:_{Q(R)}I$. Thus $JI=xI$, hence $I\subset x:J$. But $I \supset x:J$ by \Cref{traceformula}, so $(2)$ holds. 
\end{proof}

\begin{ex}
The implication $(4)\implies (1)$ does not hold even if $R$ is a hypersurface.
Let $R=k[[x,y]]/(x^5-y^3)$ and $J=(y,x^4)$. Then $(y):J=\m=\tr(J)$ but $\m \ncong \m^2$. 

\end{ex}

\begin{prop}\label{proptrI}
Suppose $I$ is a stable ideal. Then $\tr{I}\cong I^*$ and $\tr{I}=x:I$ for any principal reduction $x$ of $I$. 
\end{prop}

\begin{proof}
Let $S=I:_{Q(R)}I$. Then $I\cong S$, so $\tr{I} = \tr{S}\cong S^*\cong I^*$. The last assertion follows from \Cref{IvJ}. 
\end{proof}

\section{Stable trace ideals and reflexive ideals}\label{sec_rf}

In this section we investigate the relationship between stable trace and reflexive ideals. Continue to assume $(R,\m)$ is a local  Cohen-Macaulay ring of dimension one. An $R$-module $M$ is called self-dual if $M\cong M^*$. Evidently a self-dual module is reflexive.  We start by recalling a couple of results and definitions from \cite{dms}. 

\begin{prop}
Any regular stable trace ideal is reflexive. There is a one-to-one correspondence between the set of (regular) stable trace ideals of $R$ and the set of birational extensions $S$ of $R$ that are self-dual. 
\end{prop}

\begin{proof}
A regular stable trace ideal is self-dual by \Cref{prop_trace}. The second statement is a combination of \cite[Lemma 2.8]{dms} and \Cref{prop_trace}.
\end{proof}

\begin{dfn}
Assume $I$ is a regular (fractional) ideal of $R$. Let $b(I)$ denote the conductor of $B_R(I)$. 
\end{dfn}

\begin{dfn}
Let $M$ and $N$ be $R$-modules. We say $M$ generates $N$ is there exists a surjection \[ M^{(\Lambda)} \twoheadrightarrow  N \]
for some index set $\Lambda$. 
\end{dfn}
\begin{lem} \label{IgenJ}
Let $I$ and $J$ be ideals. Then $I$ generates $IJ$. \end{lem}

\begin{proof}  For each  element $j \in J$ we may construct a family of homomorphisms $\alpha_j: I \longrightarrow IJ$ given by $i \mapsto ij$. The product of these homomorphisms is a surjective map
\[ \prod_{j \in J} \alpha_j : \prod_{j \in J} I \longrightarrow IJ .\]
\end{proof}

\begin{rem}\label{rmkIJ}
It follows from Lemma \ref{IgenJ} that $I$ generates $I^n$ for all $n>0$. Note, when $J$ is finitely generated by elements $j_i, \ldots, j_n$ then  $\prod_{j \in J} \alpha_j$ may be replaced by  $\prod_{k=1}^{n} \alpha_{j_k}$ in the proof of Lemma \ref{IgenJ}. 
\end{rem}

\begin{prop}\label{Iulrich}
Let $I,J$ be regular (fractional) ideals. We have:
\begin{enumerate}
    \item It always holds that $b(I)\subseteq \tr I$.
    \item If $J\cong IJ$ then $\tr J \subseteq b(I)$. The converse holds if $J$ is reflexive. 
\end{enumerate}
\end{prop}

\begin{proof}
We have $b(I) = \tr B_R(I) =\tr I^n$ for some large enough $n$. But as $I$ generates  $I^n$ by \Cref{rmkIJ}, $\tr I \supseteq \tr I^n$.

For $(2)$, $IJ\cong J$ means $J$ is $I$-Ulrich in the sense of \cite{dms}, and the assertions follows from \cite[Corollary 4.11]{dms}.
\end{proof}

\begin{prop}
The following are equivalent:
\begin{enumerate}
    \item $b(I)= \tr I$.
    \item $\tr I\cong I^*$. 
\end{enumerate}
\end{prop}

\begin{proof}
Assume $(1)$, to prove $(2)$ we can make a faithfully flat extension of $R$ to assume that $I$ has a principal reduction $x$. Then as we always have $b(I)\subset x:I \subset \tr I$, see \cite[4.20]{dms}, $(1)$ forces $\tr I=x:I\cong I^*$. 

Assume $(2)$. Then as $\tr I\cong II^*$ we have $\tr I\cong I\tr I$, so $\tr I\subset b(I)$, and thus equality holds by \Cref{Iulrich}. 
\end{proof}

\begin{prop}
Assume that $I$ is a regular reflexive ideal. The following are equivalent:
\begin{enumerate}
    \item $I$ is stable. 
    \item $b(I)= \tr I$.
    \item $\tr I\cong I^*$. 
\end{enumerate}
\end{prop}

\begin{proof}
$I$ is stable is the same as $I\cong I^2$, which since $I$ is reflexive, is equivalent to $\tr I= b(I)$ by \Cref{Iulrich}.
\end{proof}

The next results identify reflexive ideals whose trace are stable. 

\begin{prop}\label{prop_ref}
Assume that $\tr I$ is stable. The following are equivalent:
\begin{enumerate}
    \item $I$ is reflexive. 
    \item $I\cong \tr I$.
    \item $I\cong I^*$. 
\end{enumerate}
In particular, if $I$ is reflexive and $\tr I$ is stable than so is $I$. 
\end{prop}

\begin{proof}
Since $\tr I$ is stable we have $\tr I\cong I^*$ by \Cref{IvJ}, so $(2) \implies (3)$ immediately. If $I$ is reflexive, then as $\tr I = b(\tr I)$,  by \Cref{Iulrich} we get that  $I\cong I\tr I\cong II^* \cong \tr I$. Thus $(1)\implies (2)$. Finally, $(3)\implies (1)$ is clear.    
\end{proof}

\begin{cor}
Let  $I$ a regular reflexive ideal. The following are equivalent:

\begin{enumerate}
    \item $\tr I$ is stable.
    \item $I$ is stable and self-dual. 
\end{enumerate}
\end{cor}

\begin{proof}
If $\tr I$ is stable, then $I$ is stable and self-dual by \Cref{prop_ref}. If $I$ is stable, then $\tr I\cong I^*\cong I$ by \Cref{proptrI}, thus $\tr I$ is stable. 
\end{proof}

\section{Two theorems on stable trace ideals}\label{sec_2thm}

The following  extends \cite[Corollary 3.2]{haydee_trace} for stable trace ideals. It plays a crucial role for classifying reflexive modules over Arf rings later. 

\begin{thm}\label{main1}
Suppose $\tr(M)=I$ is regular and stable. Then $I = \mathfrak{C}_{{\End_R(M)}/R}$, $I = \tr(M^*)$ and  $I\in \ad(M^*)$, that is,  $I$ is a summand of a direct sum of copies of $M^*$. 
\end{thm}

\begin{proof}
Let $S=Z(\End(M^*))$. Then $S = \End(I) = I^*$ as a subset of $Q(R)$; see \cite[Proposition 2.8, Corollary 3.24]{haydee_trace}. 

Since $M^*$ generates $M\otimes_R M^*$ which generates $\tr(M)$ we have $I = \tr(M) \subset \tr(M^*)$. Also, since $M^*$ is an $S$-module, $S$ generates $M^*$ and $\tr(M^*)\subset \tr(S)$. By assumption $I$ is stable and trace, that is $I \cong I^* = S$, thus $\tr(S) = \tr(I^*) = \tr I = I$. Altogether \[ I \subseteq \tr(M^*) \subseteq \tr(S) = I .\]

From this we know $M^*$ generates $I$. Because $I$ is stable, $I\cong S$, and it follows that a direct sum of copies of the $S$-module $M^*$ surjects onto $S$. Since $S$ is a torsionfree $R$-module, this surjection is also $S$-linear. It follows that $I\cong S\in \ad(M^*)$.
\end{proof}

Next, we study homomorphism modules of stable trace ideals. We start with some preparatory results.

\begin{lem}\label{lemHom}
Let $I,J$ be ideals. Then $\Hom_R(I,J)= \Hom_R(I,\tr(I)\cap J)$. 
\end{lem}

\begin{proof}
Any map from $I$ to $R$ must have its image in $\tr(I)$. If follows that any map from $I$ to $J$ must land in $\tr(I)\cap J$.  
\end{proof}

\begin{prop}\label{prophom1}
Suppose that $I,J$ are stable ideals. Then $$\tr{(IJ)}\cong \Hom_R(I,\tr{J})$$
\end{prop}

\begin{proof}
The natural surjection $I\otimes J\to IJ$, given by $i \otimes j \mapsto ij$,  is an isomorphism at minimal primes and therefore has a kernel of finite length. Applying $\Hom(-,R)$ and then using hom-tensor adjointness  and the fact that $\tr{J} \cong J^*$ (see \Cref{proptrI}) we get that \[(IJ)^* \cong \hom R {I\otimes J} R \cong \Hom_R(I,J^*) \cong \Hom_R(I,\tr{J}).\] 
Since $IJ$ is also stable, we are done. 

\end{proof}

\begin{thm}\label{propHom}
Suppose that $I,J,$ are stable trace ideals. Then $$\Hom_R(I,J)\cong\Hom_R(J,I) \cong \tr{(IJ)}$$
If $I\cap J$ is also a stable trace ideal then $$\tr{(IJ)}\cong I\cap J$$. 
\end{thm}

\begin{proof}
The first assertion follows from \Cref{prophom1}. For the second one, let $L=I\cap J$, then applying \Cref{lemHom} and \Cref{prophom1} we get:
$$\tr{(IJ)}\cong \Hom(I,J)\cong \Hom(I,L)\cong \Hom(L,I)\cong \Hom(L,L\cap I)=\Hom(L,L)\cong L $$
\end{proof}

\section{The set of regular trace ideals with the same integral closure}\label{sec_int}

In this section we explore the interaction between trace ideals and their integral closures. 

\begin{lem} \label{lem_dualpowers}
Let $J$ be an ideal and $\phi \in J^*$. Then $\phi^n$ is defined on $J^n$. In particular, for $x\in J$ \[\phi^n(x^n) = (\phi(x))^n.\] 
\end{lem}

\begin{proof}
Let $J \subseteq R$ be any ideal with $\phi \in J^*$. For $a_i\in J$ by $R$-linearlity

\[ \phi(a_1)\cdot \phi(a_2) = \phi(\phi(a_1) \cdot a_2) =\phi(\phi(a_1\cdot a_2)) = \phi^2(a_1a_2). \] 

By induction $\phi(a_1) \cdots \phi(a_n) = \phi^n(a_1 \cdots a_n)$ for all $n\geq 1$. Indeed,
\begin{align*}
  \phi(a_1) \cdots \phi(a_k) & =\phi^{k-1}(a_1 \cdots a_{k-1})\cdot \phi(a_k)\\
& =\phi^{k-1}(\phi(a_k) a_1 \cdots a_{k-1})\\
 & = \phi^{k-1}(\phi(a_1 \cdots a_{k-1}a_k))\\
 & = \phi^k(a_1 \cdots a_k).
\end{align*}

It follows that $\phi^n(x^n) = (\phi(x))^n$ for any $x \in J$ and $n \geq 1$. \qedhere


\end{proof}

\begin{thm} \label{traceclosure}
If $I$ is a trace ideal then so is $\bar{I}$. 
\end{thm}

\begin{proof}
Suppose $I$ is a trace ideal and $x\in \bar{I}$ such that \[ x^n +a_1x^{n-1} + \cdots + a_{n-1}x + a_n = 0\]

for some $n \geq 0$ and $a_i \in I^i$. Let $\psi \in (\bar{I})^*$. Note that since $x^{n-j}$ is in the domain of $\psi^{n-j}$ by Lemma \ref{lem_dualpowers} and $a_j$ is in the domain of $\psi^n \in \En R I$ for all $n \geq 0$ because $I$ is trace and so $(\psi|_I)^n \in I^* = \En R I$. For any $n \geq 1$ we have 
\begin{align*}
\psi^{n}(a_jx^{n-j}) &= \psi^j(\psi^{n-j}(a_jx^{n-j}) ) \\
& =  \psi^{j}(a_j \cdot \psi^{n-j}(x^{n-j}))\\
& = \psi^{j}(a_j) \psi^{n-j}(x^{n-j})\\
& = \psi^{j}(a_j) \cdot (\psi(x))^{n-j}
\end{align*}

where the final equality follows from Lemma \ref{lem_dualpowers}. We may write $a_j = b_1 \cdots b_j$ with $b_i \in I$ and observe that 
\[\psi^{j}(a_j) = \psi^{j}(b_1 \cdots b_j) = b_1 \cdots b_{j-1} \psi^{j}(b_j)  \in I^j.\]

By the above
\begin{align*}
(\psi(x))^n &= \psi^n(x^n)  \\
&= \psi^{n}(-(a_1x^{n-1} + \cdots + a_{n-1}x + a_n))\\
& =- \sum_{j=1}^n\psi^{j}(a_j) \cdot (\psi(x))^{n-j}
\end{align*}

It follows that $\psi(x) \in \bar{I}$ and $\bar{I}$ is a trace ideal. 
\end{proof}

\begin{lem} \label{stabletrace}

If $I$ is a regular trace ideal with a principal reduction $(x)$ and $J$ is another trace ideal with $I\subset J \subset \bar I$. If $J$ is stable, then $I=J$.
\end{lem}

\begin{proof}   
Let $(x)$ be a principal reduction of $I$. Then $(x)$ is also a principal reduction of $J$ which is stable by assumption and therefore $xJ=J^2$; see \cite[Proposition 4.5]{dms} .


It follows that $I^2\subseteq J^2\subseteq (x)$, so by Proposition \ref{prop_trace} $I^2=xI$. Recall that $J$ is also a trace ideal by  Proposition \ref{traceclosure}. Now again by \Cref{prop_trace}, $I = (x):I \supset (x):J = J$, so $I=J$.  
\end{proof}

\begin{cor} \label{corst}

If $I$ is a regular trace ideal with a principal reduction $(x)$ and $\bar I$ is stable, then $I=\bar I$.
\end{cor}

Next, we make the key definition of this section. 

\begin{dfn}
Let $J$ be a regular integrally closed ideal. Let $T(J)$ denote the set of trace ideals $I$ such that $\bar I=J$. 
\end{dfn}

To keep things simple we assume:
\begin{setup}\label{setup3}
Let $R$ be a Cohen-Macaulay local ring of dimension one with the following properties: any regular ideal has a principal reduction and the integral closure $\bar R$ of $R$ is a finite $R$-module. Let $\cond$ be the conductor of $\bar R$ in $R$. For instance, if $R$ is a analytically unramified ring with inifinite residue field then $R$ satisfies these conditions. 
\end{setup}

\begin{prop}\label{propTJ}
Let $R$ as in \Cref{setup3}. Let $J$ be a regular integrally closed ideal. We have:
\begin{enumerate}
    \item $T(J)$ is non-empty if and only if $J$ contains $\cond$. 
    \item $T(J)=\{J\}$ if $J$ is stable and contains $\cond$. 
    \item Any stable ideal in $T(J)$ is minimal with respect to inclusion.
    \item The set of regular trace ideals of $R$ is finite if and only if $T(J)$ is finite for each regular integrally closed ideal $J$. 
\end{enumerate}
\end{prop}

\begin{proof}
Any regular trace ideal contains $\cond$, and if $J\supset \cond$, then $J$ is a trace ideal, see \cite[3.5, 3.11]{dms}. Thus $(1)$ follows. For $(2)$ if  $J$ contains $\cond$ then $J\in T(J)$ and $T(J)=\{J\}$ by \Cref{corst}.   $(3)$ is just  another way to state \Cref{stabletrace}. 

Lastly, $T(J)$ is non-empty only if $J$ contains the conductor, but the set of integrally closed ideals containing the conductor is finite, \cite[6.4]{dms}, thus $(4)$ follows. 
\end{proof}

\Cref{propTJ} motivates the:

\begin{ques}\label{quesTJ}
When is $T(J)$ finite? When is $T(J)=\{J\}$?
\end{ques}

The following examples elucidate the concept of $T(J)$ and indicate why the answers to \Cref{quesTJ} might not be simple. 

\begin{ex}
Let $R=\mathbb C[[x,y]]/(xy(x-y))$ with $\m=(x,y)$. The conductor is $\m^2$. Then $T(\m)=\{\m\}$. Indeed any regular trace ideal $I$ must contain the conductor $\m^2$, and if it is not $\m$ then $I= (\m^2,l)$ with some reduction $l$ of $\m$. As $\m^3\subset (l)$, we have that $l:I=\m$. Thus the trace of $I$ contains $\m$, see \Cref{traceformula}, so $I$ is not a trace ideal.  Note that $\m$ is not stable, hence $T(J)=\{J\}$ does not necessarily imply $J$ is stable. 
\end{ex}

\begin{ex}
Let $R=\mathbb C[[x,y]]/(x^n-y^n)$. Let $l=x-ay$, where $a\in \mathbb C$ such that $a^n\neq 1$. Then $l$ is a reduction of $\m=(x,y)$.  Any ideal containing $l$ must be $I_s=(l, y^s)$ for some $s\leq n$. Clearly $l:I_s=I_{n-s}$ and $I_s$ is a trace ideal iff $s\leq n/2$. The smallest such ideal is $I_{\floor{n/2}}$. It is stable if and only if $n$ is even. Thus the minimal elements in $T(J)$ can be infinite, and they might or might not be stable. When $n\geq 4$ is even, this example also gives a ring with infinitely many stable trace ideals. 
\end{ex}

\begin{ex}
Let $R=\mathbb C[[t^4,t^5,t^6]]$ with $I_a=(t^4-at^5, t^6)$ for $a\in \mathbb C$ and $\m=(t^4,t^5,t^6)$. Then $T(\m)=\{\m\} \cup \{I_a\}_{a\in \mathbb C}$, see \cite[Example 3.4]{GotoIsobeKumashiro}. The minimal elements are $I_a$, each of them is stable. Indeed, $x=t^4-at^5$ is a minimal reduction of $I_a$. We just need to check $t^{12}\in xI$. But $t^{12}/x= \frac{t^8}{1-at}= t^8(1+at+a^2t^2+...) \in I$ because for any $s\geq 8$, $t^s\in \m^2\subset I$, as $\m^2$ is the conductor.   
\end{ex}

\section{Applications to Arf Rings} \label{SectionArf}
In \cite{KobayashiTakahashi}, Kobayashi and Takahashi describe the rings for which every ideal is isomorphic to a trace ideal. If $R$ is a commutative Noetherian ring of depth one, for example,  every ideal is isomorphic to a trace ideal if and only if $R$ is a hypersurface of Krull dimension one and multiplicity at most two; see \cite[Theorem 1.2]{KobayashiTakahashi}. Such a ring is an Arf ring; see definition below. In this section we study regular trace ideals over Arf rings in general with repercussions for candidates for test ideals for projective dimension and the structure of reflexive modules over these rings. For instance, over Arf rings with infinite residue field any reflexive ideal is isomorphic to it's own trace, and this can be viewed an extension of the aforementioned result of Kobayashi and Takahashi. We also point out several applications to rigidity of $\Tor$ and structure of reflexive modules over Arf rings.  

\begin{dfn}
Let $R$ be a 1-dimensional Cohen-Macaulay ring. $R$ is called an \emph{Arf} ring provided every integrally closed regular ideal is stable. 
\end{dfn}

\begin{ex}  
Indeed, any hypersurface ring of multiplicity less than or equal to 2 is an example of an Arf ring. Such rings are two-generated rings, over which all ideals are stable; see \cite[Theorem 3.4]{SallyVasc1}. For example, $k[[x,y]]/(y^2 - x^3)$ is an Arf Ring.
\end{ex}

\begin{ex}
Let $k$ be a field. Fix integers $e\geq 2, n\geq 1$. Let $H$ be the numerical semigroup generated by $\{e, ne+1, ne+2,..., ne+{e-1}\}$. Then $R_{e,n} = k[[t^a, a\in H]]$ is an Arf local domain with multiplicity $e$ and embedding dimension $n$ (see \cite[Example 4.2]{daoArf}). 
\end{ex}

\begin{thm}\label{mainT1}
Let $R$ be a Cohen-Macaulay local ring of dimension one such that any regular ideal has a principal reduction. The following are equivalent:
\begin{enumerate}
    \item $R$ is Arf.
    \item Any regular trace ideal is stable.  
\end{enumerate}
\end{thm}

\begin{proof}
Assume $R$ is Arf, and let $I$ be a regular trace ideal. Then $\bar I$ is stable, which forces $I=\bar I$ by \Cref{corst}. Thus $I$ is stable. 

Assume $(2)$. Let $J$ be a regular integrally closed ideal and let $I=\tr(J)$. Then $I$ is stable, so $I\cong J^*$ by \Cref{IvJ}. Then as $J$ is reflexive (\cite[3.11]{dms}), $J\cong I^* \cong I$ by \Cref{prop_trace}. As $I$ is stable, so is $J$. Thus $R$ is Arf by definition. 

\end{proof}

\begin{cor}
Let $R$ be an Arf ring such that any regular ideal has a principal reduction. Then any regular reflexive ideal is isomorphic to a trace ideal (necessarily it's own trace). 
\end{cor}

\begin{proof}
The result follows from \Cref{prop_ref} and \Cref{mainT1}.
\end{proof}

\begin{prop} \label{traceintegrallyclosed}
Let $R$ be an Arf ring such that any regular ideal has a principal reduction. Then the set of regular trace ideals and integrally closed ideals containing the conductor coincide. If $\overline R$ is finite then the set of regular trace ideals is finite. 
\end{prop}

\begin{proof}
If an ideal is integrally closed and contains the conductor it is a regular trace ideal; see \cite[3.11]{dms}. By assumption, each regular trace ideal $I$ has a principle reduction and since $R$ is Arf ring, $\bar{I}$ is stable.  By Lemma \ref{stabletrace}, $I=\bar{I}$.  

\end{proof}

\begin{rem}
The following gives an alternative proof to the key results in \cite{daoArf} and \cite{isobe2021reflexive}.
\end{rem}

\begin{prop} \label{refltrace}
Let $R$ be a complete Arf ring such that any regular ideal has a principal reduction. 
Then any reflexive module is isomorphic to a direct sum of integrally closed ideals. 
\end{prop}

\begin{proof}
Any reflexive module $M$ is the direct sum of indecomposable reflexive summands, each of which is the $R$-dual of some other module. By Proposition 3.5 and Krull-Remak-Schmidt, $M$ is therefore a direct sum of trace ideals. Because $R$ is Arf, Lemma 3.6 implies these trace ideal summands are also stable and integrally closed.
\end{proof}

\begin{cor}
Let $R$ be an Arf ring such that any regular ideal has a principal reduction. Then  if $M$ is reflexive and $\Tor^R_t(M, N ) = 0$ then $pd_R(N) \leq t$.  
\end{cor}

\begin{proof}
 The result follows directly from Corollary \ref{refltrace} and  \cite[Theorem 2.10]{idealcaseHWC}; see also \cite[Corollary 3.3]{corso2005integral}. 
\end{proof}






\begin{cor}
For complete Arf domains there are only finitely many classes of indecomposable reflexive modules.

\end{cor}
\begin{proof}
The result follows from Proposition \ref{refltrace} and \cite[6.4]{dms}.
\end{proof}

\begin{cor}
Let $R$ be an Arf ring and $I,J$ are integrally closed ideals containing the conductor. Then $$\Hom_R(I,J)\cong\Hom_R(J,I) \cong (IJ)^* \cong \tr{(IJ)} \cong I\cap J.$$
\end{cor}

\begin{proof}
By Proposition \ref{traceintegrallyclosed}, the integrally closed ideals containing the conductors are precisely the regular trace ideals of $R$. Since the intersection of integrally closed ideals is also integrally closed, the ideals $I$, $J$ and $I \cap J$ are all stable trace ideals. Thus we are done by \Cref{propHom}.

\end{proof}

\section{Applications to rings defined by trace of the canonical module}\label{sec_gor}

Recently, there have been a lot of activities around classes of rings defined by trace of the canonical module, see \cite{herzog2019trace, herzog2021tiny, dao2020trace}. We point out that many of the properties of such rings can be recovered by our results. For instance,  we recover a result in \cite{herzog2019trace} showing when the classes of almost Gorenstein rings and nearly Gorenstein rings coincide. We also prove a generalization of the fact in \cite{herzog2021tiny} that the trace of a canonical ideal $\omega$ is equal to the conductor if and only if $\omega^2$ is isomorphic to the conductor.

\begin{cor} \label{containment}
Suppose $R$ has minimal multiplicity and $J$ is a regular ideal with principal reduction $x$. Then $(x):J \supset \m$ if and only if $\tr(J)\supset \m$. 
\end{cor}

\begin{proof}
Assume $(x): J \supseteq \mathfrak{m}$. Since $x \in J$, note that \[\mathfrak{m} \subseteq (x):J = (x) J^* \subseteq JJ^* =  \tr J.\] 

On the other hand, if $\tr(J) \supseteq \mathfrak{m}$ then $\tr (J) \in \{ \mathfrak{m}, R \}$. In either case $\tr(J) = \tr(J)^2$ because either $R = R^2$ or, since $R$ has minimal multiplicity, $\mathfrak{m} = \tr(J)$ is a stable and as a stable trace ideal  $\mathfrak{m} \cong \mathfrak{m}^2$; see by Proposition \ref{prop_trace}. By Proposition \ref{IvJ} it follows that $ (x): J = \tr (J) \supseteq \mathfrak{m}$. 
\end{proof}

\begin{dfn}
Let $R$ be a local Cohen-Macaulay rings with maximal ideal $\mathfrak{m}$ and  canonical module $\omega_R$. We say $R$ is \emph{almost Gorenstein} if \[ (x): \omega_R \supseteq \mathfrak{m} ,\]

for some principal reduction $(x)$ of $\omega_R$.

We say $R$ is \emph{nearly Gorenstein} provided $\tr{\omega_R} \supset \mathfrak{m}.$
\end{dfn}

\begin{rem}
If $J$ is a canonical ideal, Corollary \ref{containment} says that being almost Gorenstein and nearly Gorenstein are equivalent. Thus we recover \cite[Theorem 6.6]{herzog2019trace}:
\end{rem}

\begin{cor}
Let $(R, \mathfrak{m})$ be a Cohen–Macaulay local ring of minimal multiplicity with $\dim R = 1$
and infinite residue field. Assume that R possesses the canonical module $\omega_R$.
If R is nearly Gorenstein, then it is almost
Gorenstein.
\end{cor}

\begin{proof} 
Since $R$ has an infinite residue field, each ideal has a principal reduction; see \cite[ 8.3.7, 5.1.6]{SwansonHuneke}. 
The result follows by applying Corollary \ref{containment} in the case $J = \omega_R$.
\end{proof}

Finally, we study ideals whose trace is the conductor. 

\begin{lem}\label{lemcantrace}
Assume \Cref{setup3}. Let $\omega$ be a canonical ideal of $R$ and $I$ a regular ideal. The following are equivalent:
\begin{enumerate}
    \item $\tr{I}=\cond$.
    \item $I^*\cong \cond$.
    \item $I\omega\cong \cond$.
\end{enumerate}
\end{lem}

\begin{proof}
Note that $\cond$ is a stable trace ideal, so $(1)\implies (2)$ by \Cref{IvJ}. Assume $(2)$, then $\tr(I) \cong II^*\cong I\cond \cong \cond$ (see \Cref{Iulrich}), so $\tr(I)=\cond$ as they are both trace ideals. 

The map $I\otimes \omega \to I\omega$ is surjective with finite length kernel (as $I$ is regular and hence locally free on the minimal primes). Take $\Hom_R(-,\omega)$ and use Hom-tensor adjointness we get that $$\Hom_R(I\omega, \omega) \cong I^*$$. 

So $(2)$ is equivalent to $\Hom_R(I\omega, \omega)\cong \cond$, which is equivalent to $I\omega\cong \Hom_R(\cond, \omega)$. But as: $$\Hom_R(\cond, \omega)\cong \Hom_R(\bar R, \omega) \cong \bar R\cong \cond $$
we are done.

\end{proof}

The following is immediate from \Cref{lemcantrace}, and recover part of \cite[Theorem 2.5]{herzog2021tiny}.

\begin{cor}
Assume \Cref{setup3}. Let $\omega$ be a canonical ideal of $R$. Then $\tr(w)=\cond$ if and only if $\omega^2\cong \cond$. 
\end{cor}

\bibliographystyle{siamplain}	
\bibliography{myrefs}		


\end{document}